\documentclass[3p,sort&compress]{elsarticle}
\usepackage{amsmath,amssymb,amsthm}
\usepackage{graphicx}
\usepackage{verbatim}

\journal{Journal of Computational Physics}
\begin{document}

\begin{frontmatter}
\title{Discretization of Fractional Differential Equations by a Piecewise Constant Approximation}

\author[UNSW]{C. N. Angstmann}
\author[UNSW]{B. I. Henry\corref{cor1}}
\ead{b.henry@unsw.edu.au}
\author[WITS,COE]{B. A. Jacobs}
\author[UNSW]{A. V. McGann}

\address[UNSW]{School of Mathematics and Statistics, UNSW Australia, Sydney NSW 2052 Australia}
\cortext[cor1]{Corresponding author}
\address[WITS]{
	School of Computer Science and Applied Mathematics, 
			University of the Witswatersrand, 
			Johannesburg, 
			Private Bag 3, 
			Wits 2050, 
			South Africa}
\address[COE]{
			DST-NRF Centre of Excellence in Mathematical and Statistical Sciences (CoE-MaSS),
            University of the Witwatersrand, 
            Johannesburg,
            Private Bag 3, 
            Wits 2050, 
            South Africa}

\date{\today}

\begin{abstract}
There has recently been considerable interest in using a nonstandard piecewise approximation to formulate fractional order differential equations as difference equations that describe the same dynamical behaviour and are more amenable to a dynamical systems analysis. Unfortunately, due to mistakes in the fundamental papers, the difference equations formulated through this process do not capture the dynamics of the fractional order equations. We show that the correct application of this nonstandard piecewise approximation leads to a one parameter family of fractional order differential equations that converges to the original equation as the parameter tends to zero. A closed formed solution exists for each member of this family and leads to the formulation of a difference equation that is of increasing order as time steps are taken. Whilst this does not lead to a simplified dynamical analysis it does lead to a numerical method for solving the fractional order differential equation. The method is shown to be equivalent to a quadrature based method, despite the fact that it has not been derived from a quadrature. The method can be implemented with non-uniform time steps. An example is provided showing that the difference equation can correctly capture the dynamics of the underlying fractional differential equation. 
\end{abstract}
\end{frontmatter}

\section{Introduction}

There has been an increasing use of fractional calculus in models of physical systems, such as anomalous diffusion \cite{MK2000, MS2012,ADH2013mmnp}, viscoelasticity \cite{M2010},and the spread of disease \cite{AHM2016fiSIR,AHM2016frSIR}. These models are often posed in the form of fractional differential equations (FDEs). Recently there have been attempts to construct a difference equation that captures the dynamical behaviour of a FDE by using a nonstandard piecewise approximation \cite{el2013discretization,Agarwal2013,ismail2015generalized,selvam2016dynamics,el2014discretization,EES2014discretization,el2016fractional}. These attempts have resulted in the incorrect  construction of first order difference equations that can not approximate the dynamics of FDEs. In the earliest of these works \cite{el2013discretization} the first order difference equation is the result of the use of incorrect limits on an integral, (see Section 2 of \cite{el2013discretization}). Later papers, such as \cite{Agarwal2013}, correct this mistake but introduce other, more subtle, errors that result in the same incorrect first order difference equations.

Traditional discretization methods for fractional differential operators, such as the Gr\"unwald-Letnikov approximation and its generalizations \cite{P1999,BDST2012,BKS2015,Jacobs2015}, result in sums over the entire history of the function. The memory effect of fractional order equations has been physically corroborated and is imperative that a numerical scheme capture this dynamic, as in the physically consistent numerical schemes derived in \cite{ADHN2015,angstmann2016stochastic}. This history dependent dynamics cannot be captured by a first order difference equation. This alone is enough to question the validity of a first order difference equation from a discretization scheme for a fractional differential equation. 

Here we consider an initial value FDE of the form,
\begin{equation}
\label{eq_f}
	_{\mathrm{C}} \mathcal{D}_{0,t}^{\alpha} x(t) = f(x(t)),\; x(0)=x_0,
\end{equation}
where $_{\mathrm{C}} \mathcal{D}_{0,t}^{\alpha}$ is a Caputo fractional derivative with $0< \alpha<1$ and $f(x(t))$ is potentially a nonlinear vector field. We show that the first order discretization given by El-Sayed and others \cite{el2013discretization,Agarwal2013,ismail2015generalized,selvam2016dynamics,el2014discretization,EES2014discretization,el2016fractional} is incorrect. Further to this we correct the derivation of the difference equation resulting in an increasing order difference equation which captures the memory effect of the FDE. This is achieved by a piecewise constant approximation of the vector field. Resulting in a one parameter family of integrable FDEs that limit to the original FDE. The integrable FDEs have a closed form solution that can be discretised to provide a difference equation that approximates the solution of the original FDE. Furthermore, we show that this method may be implemented with an non-uniform time step. An example is presented that shows the difference equation correctly captures the dynamics of a specific FDE.

\section{Fractional Derivatives}
The properties of fractional derivatives have been widely explored in \cite{OS1974,P1999} and more recently in \cite{LD2007}. There exist multiple types of fractional derivatives, here we will focus on FDEs involving Caputo derivatives. A Caputo fractional derivative is defined by \cite{C1967},
\begin{equation}
_{\mathrm{C}}\mathcal{D}^{\alpha}_{0,t}x(t)=\frac{1}{\Gamma(1-\alpha)}\int_{0}^{t}(t-t')^{-\alpha}\frac{d x(t')}{d t'}dt' ,
\end{equation}
for $0<\alpha<1$. 
The Riemann-Liouiville fractional derivative is defined by,
\begin{equation}
_{\mathrm{RL}}\mathcal{D}^{\alpha}_{0,t}x(t)=\frac{1}{\Gamma(1-\alpha)}\frac{d}{dt}\int_{0}^{t}(t-t')^{-\alpha}x(t')dt', 
\end{equation}
for $0<\alpha<1$. 
We can transform a Caputo derivative to a Riemann-Liouville through the following relation \cite{LD2007};
\begin{equation}\label{eq:RL_C}
_{\mathrm{RL}}\mathcal{D}^{\alpha}_{0,t}(x(t)-x(0))= \, _{\mathrm{C}}\mathcal{D}^{\alpha}_{0,t}x(t),
\end{equation}
hence the two are equivalent when $x(0)=0$.

Fractional derivatives can be easily expressed in Laplace space and we will make use of this form.

The Laplace transform of the Caputo fractional derivative is,
\begin{equation}
\mathcal{L}_t\{_C\mathcal{D}^\alpha_{0,t}x(t)\}=s^\alpha\mathcal{L}_t\{x(t)\}-s^{\alpha-1}x(0).
\end{equation}

The Gr\"unwald-Letnikov derivative is defined by,
\begin{equation}
_{\mathrm{GL}}\mathcal{D}^{\alpha}_{0,t}x(t)=\lim_{\delta t\to0}\frac{1}{\delta t^{\alpha}}\sum_{m=0}^{\infty}(-1)^m \left( \begin{array}{c} \alpha \\  m \end{array} \right)x(t-m \delta t).
\end{equation}
If $x(t)\in C^0$ and $0<\alpha\leq1$, then this is equivalent to the Riemann-Liouville derivative. 
Gr\"unwald-Letnikov derivatives have long been used as a basis of methods for discretising Riemann-Liouville FDEs \cite{P1999,BDST2012}.

\section{Comparison to El-Sayed and Salman}

In general the approach taken by El-Sayed and others \cite{el2013discretization,Agarwal2013,ismail2015generalized,selvam2016dynamics,el2014discretization,EES2014discretization,el2016fractional} can be generalised easily for an equation in the form of Eq. (\ref{eq_f}). This gives a discrete equation of the form,
\begin{equation}
\label{eq_wrong}
x((n+1)\delta t)=x(n\delta t)+\frac{\delta t^{\alpha}}{\Gamma(1+\alpha)}f(x(n\delta t)).
\end{equation} 
This first order difference equation can not capture the dynamics of a fractional order differential equation.

\subsection{Convergence}
It is easy to see that Eq. (\ref{eq_wrong}) does not have a well defined limit as $\delta t\to0$ such that $t=n \delta t$. We can rearrange Eq. (\ref{eq_wrong}) to give,
\begin{equation}
\Gamma(1+\alpha)\frac{(x(t+\delta t)- x(t))}{\delta t^\alpha}=f(x(t)).
\end{equation} 
Taking $x(t+\delta t)=x(t)+\delta t x'(t)+\frac{\delta t^2}{2}x''(t)+o(\delta t^2)$, this becomes,
\begin{equation}
\Gamma(1+\alpha)\frac{(\delta t x'(t)+\frac{\delta t^2}{2}x''(t)+o(\delta t^2))}{\delta t^\alpha}=f(x(t)).
\end{equation} 
Finally taking the limit $\delta t \to 0$ gives $f(x(t))=0$,
for $0<\alpha<1$. Hence we do not recover the original equation in the limit of small time steps. It should also be noted that the correct equation is recovered in the case $\alpha=1$.

\section{Piecewise Constant Integrablization}

Whilst the final result of the El-Sayed method \cite{el2013discretization,Agarwal2013,ismail2015generalized,selvam2016dynamics,el2014discretization,EES2014discretization,el2016fractional}  is incorrect, the initial approach that they undertake has some merit. Here we will consider the general Caputo initial value problem given in Eq. (\ref{eq_f}). We will construct a family of equations, parametrized by $\delta t$ such that in the limit $\delta t \rightarrow 0$ the family limits to Eq. (\ref{eq_f}). Each member of the family is an integrable FDE with a closed form solution. We will refer to this process as an integrablization of Eq. (\ref{eq_f}). Here this is achieved by replacing the right hand side of the FDE with a piecewise constant function.  Choosing some time step $\delta t$, the FDE to solve would become,
\begin{equation}
\label{eq_d}
_{\mathrm{C}}\mathcal{D}^{\alpha}_{0,t}x(t)=f\left(x\left(\delta t \left\lfloor\frac{t}{\delta t}\right\rfloor\right)\right).
\end{equation}
This is the same equation that El-Sayed and Salman \cite{el2013discretization} attempt to solve. We show below that this is an integrable equation whose solution is trivially obtained. 
The piecewise constant function is chosen so that in the limit of small $\delta t$ we recover the original equation, i.e.
\begin{equation}
\lim_{\delta t \to 0} f\left(x\left(\delta t \left\lfloor\frac{t}{\delta t}\right\rfloor\right)\right)= f(x(t)).
\end{equation}
Using a unit step function, defined by
\begin{equation}
u(t) = 
\left\{
\begin{array}{ll}
      0 & t < 0, \\
      1 & t \geq 0, \\
\end{array} 
\right.
\end{equation}
we can rewrite the right-hand side Eq. (\ref{eq_d}) as a sum, giving 
\begin{equation}
\label{eq_piececonst}
_{\mathrm{C}}\mathcal{D}^{\alpha}_{0,t}x(t)=\sum_{m=0}^{\infty}f(x(m \delta t))(u(t-m\delta t)-u(t-(m+1)\delta t)).
\end{equation}
This can be expressed as a sum of single unit step functions,
\begin{equation}
\label{eq_piececonstre}
_{\mathrm{C}}\mathcal{D}^{\alpha}_{0,t}x(t)=f(x_0)+\sum_{m=1}^{\infty}\left(f(x(m \delta t))-f(x((m-1)\delta t))\right)u(t-m\delta t).
\end{equation}
The solution of this equation can be found utilising Laplace transforms. The Laplace transform of Eq. (\ref{eq_piececonstre}) yields,
\begin{equation}
s^{\alpha}\mathcal{L}\{x(t)\}-s^{\alpha-1}x_0 = s^{-1}f(x_0)+s^{-1}\sum_{m=1}^{\infty}\left(f(x(m \delta t))-f(x((m-1)\delta t))\right)e^{-s m \delta t},
\end{equation}
where $x(0)=x_0$. Rearranging and inverting the Laplace transform gives,
\begin{equation}
x(t)=x_0+\frac{t^{\alpha}}{\Gamma(1+\alpha)}f(x_0)+\sum_{m=1}^{\infty}\left(\frac{(t-m\delta t)^{\alpha}}{\Gamma(1+\alpha)}(f(x(m\delta t))-f(x((m-1)\delta t)))\right)u(t-m\delta t).
\end{equation}
We note that this is a solution in continuous $t$. This can be simplified to an $n^{\mathrm{th}}$ order difference equation by setting $t=n\delta t$,  
\begin{equation}
\label{eq_pwc}
x(n \delta t)=x_0+\frac{(n\delta t)^{\alpha}}{\Gamma(1+\alpha)}f(x_0)+\sum_{m=1}^{n-1}\frac{((n-m)\delta t)^{\alpha}}{\Gamma(1+\alpha)}(f(x(m\delta t))-f(x((m-1)\delta t))).
\end{equation}
This difference equation incorporates the history of the function and as such the order of the difference equation grows with each time step.

\section{Convergence}

It is easy to confirm that the discretization given in Eq. (\ref{eq_pwc}) will limit to the solution of Eq. (\ref{eq_f}) in the limit as $\delta t\to 0$. Firstly we note that the solution of Eq. (\ref{eq_f}) can be found by fractionally integrating both side of the equation giving,
\begin{equation}
x(t)-x(0)=\int_0^t\frac{(t-t')^{\alpha-1}}{\Gamma(\alpha)}f(x(t'))dt'.
\end{equation}
Integrating by parts leads to,
\begin{equation}
\label{eq_sol}
x(t)=x(0)+\frac{t^\alpha f(x(0))}{\Gamma(1+\alpha)}+\int_0^{t}\frac{(t-t')^\alpha}{\Gamma(1+\alpha)}\frac{d}{dt'}f(x(t'))dt'.
\end{equation}
Now considering Eq.(\ref{eq_pwc}), we can write,
\begin{equation}
x(n \delta t)=x_0+\frac{(n\delta t)^{\alpha}}{\Gamma(1+\alpha)}f(x_0)+\sum_{m=1}^{n-1}\frac{\delta t((n-m)\delta t)^{\alpha}}{\Gamma(1+\alpha)}\frac{(f(x(m\delta t))-f(x((m-1)\delta t)))}{\delta t}.
\end{equation}
Taking the limit $\delta t \to 0$ such that $t=n \delta t$ and $t'=m\delta t$ are fixed then one recovers Eq. (\ref{eq_sol}). This shows that the discretization given in Eq. (\ref{eq_pwc}) converges to the solution of Eq. (\ref{eq_f}) in the limit $\delta t \to 0$. This also shows that the discretization could have been derived from a quadrature of Eq. (\ref{eq_sol}). This integral form representation of the method shows that this is related to a fractional order Adams method \cite{ZM2016}.

\section{Non-uniform Time Step}

One benefit of this approach is the ease of implementing a non-uniform time step. The approximation to Eq. (\ref{eq_f}), can be formulated so that the time steps, $\{\delta\tau_1, \delta\tau_2\, \ldots \}$, are not uniformly sized. Let the sum of the first $i$ time steps be represented as $\tau_i$ i.e. 
\begin{equation}
\tau_i = \sum_{j=1}^{i} \delta \tau_j.
\end{equation}
With the unit step notation the FDE approximation is,
\begin{equation}
\label{eq_piececonstNon}
_{\mathrm{C}}\mathcal{D}^{\alpha}_{0,t}x(t)=f(x_0)+\sum_{m=1}^{\infty}\left(f(x(\tau_m))-f(x(\tau_{m-1}))\right)u(t-\tau_m).
\end{equation}
Again, the solution can be found using Laplace transform techniques. Following the same method as above this gives, 
\begin{equation}
x(t)=x_0+\frac{t^{\alpha}}{\Gamma(1+\alpha)}f(x_0)+\sum_{m=1}^{\infty}\left(\frac{(t-\tau_m)^{\alpha}}{\Gamma(1+\alpha)}(f(x(\tau_m))-f(x(\tau_{m-1})))\right)u(t-\tau_{m}).
\end{equation}
This leads to a difference equation by setting $t=\tau_n$,  
\begin{equation}
\label{eq_pwcNon}
x(\tau_n)=x_0+\frac{\tau_n^{\alpha}}{\Gamma(1+\alpha)}f(x_0)+\sum_{m=1}^{n-1}\frac{(\tau_n-\tau_m)^{\alpha}}{\Gamma(1+\alpha)}(f(x(\tau_m)-f(x((\tau_{m-1}))).
\end{equation}
This difference equation reduces to Eq. (\ref{eq_pwc}) when uniform time steps are taken.

\section{Examples}

\subsection{Comparision to El-Sayed et. al.}

In \cite{el2013discretization} El-Sayed and Salman consider the fractional order Riccati equation of the form,
\begin{equation}
\label{eq_Riccati}
_{\mathrm{C}}D^{\alpha}_{0,t}x(t)=1-\rho (x(t))^2,
\end{equation}
with the initial condition $x(0)=x_0$. The first order difference equation that they derive is then given as,
\begin{equation}
\label{eq_wrongRiccati}
x((n+1)\delta t)=x(n\delta t)+\frac{\delta t^\alpha}{\Gamma(1+\alpha)}(1-\rho (x(n\delta t))^2).
\end{equation}

\begin{figure}[h]
\begin{center}
\includegraphics[width=100mm]{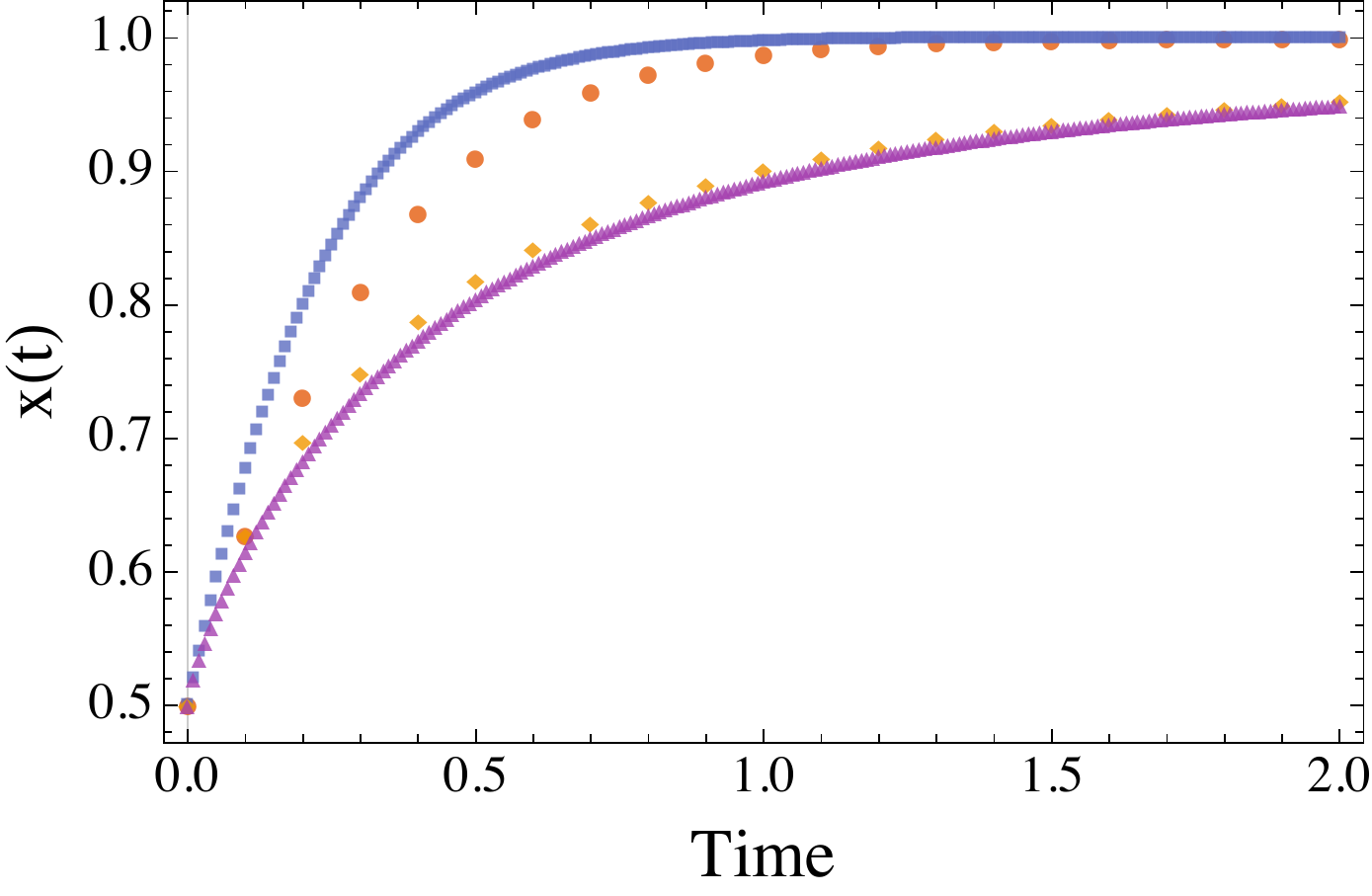}
\caption{The results of the El-Sayed discretization of the Riccati equation, Eq. (\ref{eq_wrongRiccati}) for $\Delta t=0.1$ (Orange Circles) and $\Delta t=0.01$ (Blue Squares), and the results of the integrablization, Eq. (\ref{eq_disRiccati}), for $\Delta t=0.1$ (Yellow Diamond) and $\Delta t=0.01$ (Purple Triangles).}
\label{fig_Riccati}
\end{center}
\end{figure}

Our discretization of Eq. (\ref{eq_Riccati}) can be found from Eq. (\ref{eq_pwcNon}),
\begin{equation}
\label{eq_disRiccati}
x(\tau_n)=x_0-\frac{\tau_n^{\alpha}}{\Gamma(1+\alpha)}(1-\rho x_0^2)-\sum_{m=1}^{n-1}\frac{(\tau_n-\tau_m)^{\alpha}}{\Gamma(1+\alpha)}[(1-\rho (x(\tau_m))^2)-(1-\rho (x(\tau_m-1))^2)].
\end{equation}
Taking $\tau_n=n\delta t$ leads to the fixed time step discretization. Figure \ref{fig_Riccati} shows the results of these two discritization methods with $\rho=1$, $\alpha=0.8$, and $x_0=0.5$, for $\delta t= 0.1,\;0.01,$ and $0.001$. From the figure it is not clear that Eq. (\ref{eq_wrongRiccati}) is convergent as $\delta t \to 0$, this is in contrast to the results from Eq. (\ref{eq_disRiccati}). We also see that the results from Eq. (\ref{eq_disRiccati}) display a much slower approach to the equilibrium at $x=1$ that is characteristic of fractional order differential equations.  

\subsection{Linear Fractional Differential Equation}
To show that the discretization presented here correctly captures the dynamics of a fractional order differential equation we will consider a simple case with a known exact solution. 
Consider the Caputo FDE,
\begin{equation}
\label{eq_ex1}
_{\mathrm{C}}D^{\alpha}_{0,t}x(t)=-c x(t),
\end{equation}
with $x(0) = x_0$. As this is a linear equation the solution is easily found by Laplace transform methods,
\begin{equation}
x(t)=x_0 E_{\alpha}(-c{t^\alpha}),
\end{equation}
where $E_{\alpha}(y)$ is a Mittag-Leffler function. 

To check that the discretization correctly captures the dynamics we will compare the exact solution with discrete points generated by the piecewise constant integrablization. We will also compare against the standard  Gr\"unwald-Letnikov discretization of the same equation. The Gr\"unwald-Letnikov discretization of Eq. (\ref{eq_ex1}) is given by,
\begin{equation}
x(n\delta t)= x_0-c \delta t^{\alpha} x((n-1)\delta t)-\sum_{m=1}^{n-1}(-1)^m \left( \begin{array}{c} \alpha \\  m \end{array} \right)\left(x((n-m) \delta t)-x_0\right).
\end{equation}
The difference equation from a  piecewise constant integrablization can be found from Eq. (\ref{eq_pwcNon}),  
\begin{equation}
x(\tau_n)=x_0-\frac{c(\tau_n)^{\alpha}}{\Gamma(1+\alpha)}x_0-\sum_{m=1}^{n-1}\frac{c(\tau_n-\tau_m)^{\alpha}}{\Gamma(1+\alpha)}(x(\tau_m)-x(\tau_{m-1})).
\end{equation}
Taking $\tau_n=n\delta t$ leads to the fixed time step discretization. As well as a fixed time step we will consider two cases of non-uniform time steps. In the first case we will draw a set of $\tau_n$'s from a uniform distribution such that the expected value of $\tau_n-\tau_{n-1}$ is $\delta t$. In the second case we deterministically chose the $\tau_n$ such that the difference between subsequent $\tau$'s is an increasing function. Again the time steps are chosen so that the average value is $\delta t$. 

\begin{figure}[h]
\begin{center}
\includegraphics[width=150mm]{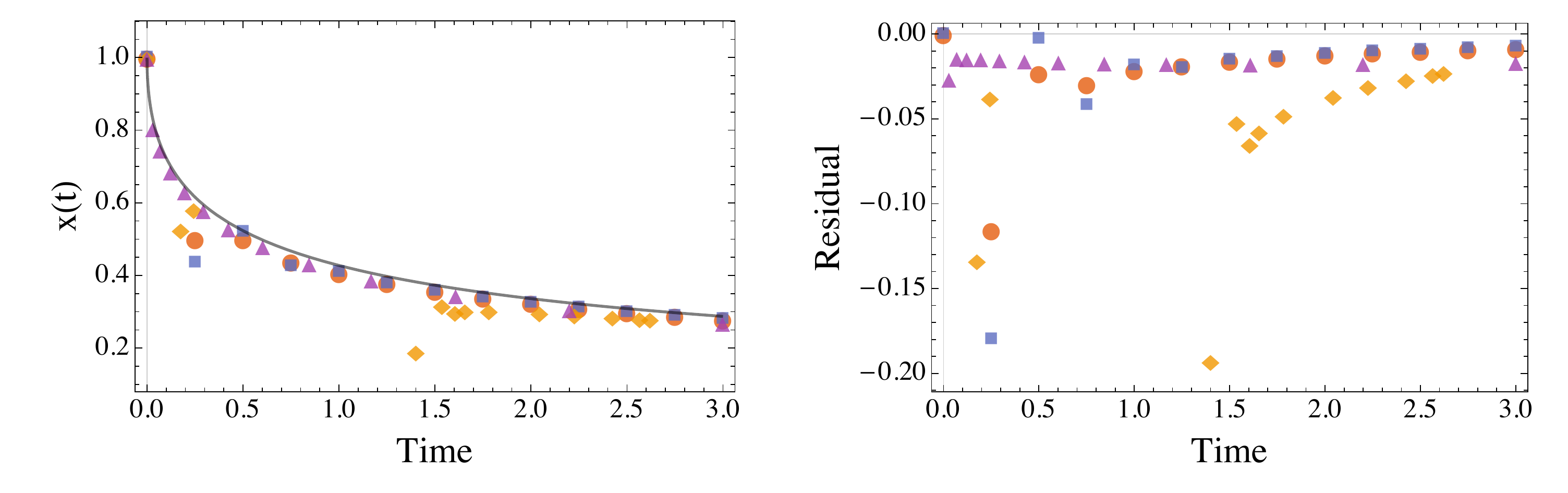}
\caption{The Gr\"unwald-Letnikov discretization (Orange Circles) and the piece wise constant integrablization on the time interval $[0,3]$ (Blue Squares, fixed time step, Yellow Diamonds, random time steps, and Purple Triangles, non-uniform spaced time steps) of Eq. (\ref{eq_ex1}) with $\alpha=0.5$, $c=1$, $x_0 =1$, and $\delta t=0.25$. The left panel shows the solutions, with the exact solution given as a solid black line, and the right panel shows the residuals.}
\label{fig1}
\end{center}
\end{figure}

Figure \ref{fig1} shows the results of these discretizations on the time interval $[0,3]$, where we have taken, $\alpha=0.5$, $c=1$, $x_0 =1$, and $\delta t=0.25$. The residuals are calculated by taking the difference between the discretization value and the exact value, i.e. $x(\tau_n)-E_{\alpha}(-c{\tau_n^\alpha})$. For a fixed time step the Gr\"unwald-Letnikov and PWC discretizations are similar in their accuracy. As we would expect taking a set of $\tau_n$'s that sample the dynamics more closely at earlier times, when the solution has a larger gradient, gives a better approximation for small times. The random sampling  time step shows that the discretization scheme is robust to the choice of time step.  

\section{Conclusion}
We have shown that the discretization method, presented by El-Sayed and others, for a class of initial values problems involving Caputo derivatives results in an incorrect first order difference equation. This difference equation, being of first order, cannot capture the dynamics of the original FDE. We have presented a correct derivation of an increasing order difference equation based on a piecewise constant approximation for the vector field of the FDE. This discretization method is amenable to non-uniform time steps and can easily be implemented on nonlinear FDEs, including non-autonomous FDEs.

\section*{Acknowledgments}
This work was supported by the Australian Commonwealth Government (ARC DP140101193). B.A.J. acknowledges the National Research Foundation of South Africa under grant number 94005. C.N.A. and B.A.J. thank the DST-NRF Centre of Excellence in Mathematical and Statistical Sciences (CoE-MaSS) for their support of C.N.A.'s visit to the University of the Witwatersrand. Opinions expressed and conclusions arrived at are those of the author and are not necessarily to be attributed to the CoE-MaSS. 

\bibliographystyle{elsarticle-num}

%\bibliography{DiscretisationFDE} 

\end{document}